\documentclass[fleqn,12pt,twoside]{article}

\usepackage{graphicx}

\title{On the Geometry of Static Space--Times}

\author{Miguel S\'anchez\address[MCSD]{Departamento de Geometr\'{\i}a y Topolog\'{\i}a, Universidad de Granada,\\
Fac. Ciencias, Avda. Fuentenueva s/n, E-18071 Granada, Spain.}
\thanks{Partially supported by a MCyT-FEDER Grant BFM 2001-2871-C04-01}}

\title{On the Geometry of Static Space--Times}
\author{M. S\'anchez}

\usepackage{amssymb} 

\newcommand{\cvd}{{\rule{0.5em}{0.5em}}\smallskip}

\begin{document}

\newtheorem{theo}{Theorem}[section]
\newtheorem{teor}[theo]{Theorem}
\newtheorem{rema}[theo]{Remark}
\newtheorem{remas}[theo]{Remarks}
\newtheorem{lema}[theo]{Lemma}
\newtheorem{propo}[theo]{Proposition}
\newtheorem{coro}[theo]{Corollary}
\newtheorem{defi}[theo]{Definition}

\font\ddpp=msbm10 scaled \magstep 1 
\newcommand\R{\hbox{\ddpp R}}    
\newcommand\N{\hbox{\ddpp N}}    
\newcommand\Lset{\hbox{\ddpp L}} 
\newcommand{\be}{\begin{equation}}
\newcommand{\ee}{\end{equation}}
\newcommand{\wk}{\rightharpoonup} 
\newcommand{\acca}{{\emph{H}}}
\newcommand{\elle}{{\emph{L}}}
\newcommand{\m}{{\emph{M}}}
\newcommand{\mo}{\m_0}   
\newcommand{\dimo}{{\bf Proof.\ }}
\newcommand{\cat}{{\rm cat}}  
\newcommand{\QED}{\vrule height .9ex width .8ex depth -.1ex}



\maketitle

\begin{abstract}
\noindent We review  geometrical properties of a static spacetime $(M,g)$, including geodesic completeness, causality, standard splittings, compact $M$, closed geodesics and geodesic connectedness. We pay special attention to the critical quadratic behavior at infinity of the coefficients $\beta$, $\beta^{-1}$  ($\beta = -g(K,K)$, being $K$ a timelike irrotational Killing vector field), which essentially control completeness, causality and geodesic connectedness. Recent references are specially discussed.  

\smallskip

\end{abstract}

%

\section{Introduction}

\noindent Static spacetimes 
are one of the simplest classes of Lorentzian manifolds. They  admit well-known physical interpretations related to the synchronizability and invariance of the metric for a field of priviledged observers $K$ (see, for example, \cite{SW}), and 
 include some classical spacetimes, as outer Schwarzschild and Reissner Nordstr\"om. Thus, many of their geometric properties have been studied from different viewpoints and, recently, there has been renewed progress made \cite{AU}, \cite{BCFS}, \cite{CMP},   \cite{Sa-St}. 
Our purpose is to review and discuss some of the geometrical properties and techniques, sketching some of the proofs (at least when the results have not been explicitly stated in the references).

This paper is organised as follows. In Section \ref{s2} both, geodesic completeness  and the possibility to split the static spacetime as a standard one, 
are characterized  (Theorems \ref{t2a}, \ref{t2b}). In Section \ref{s3},  the causal ladder of a standard static spacetime is shown to be  causally continuous but not  causally simple (Theorem \ref{t3a}, Remark \ref{r3a}), and possible standard splittings in the globally hyperbolic case are discussed, Remark \ref{r3b}. Causality Theory and different geometrical tools apply for the problem of closed geodesics and connectedness by causal geodesics. In Section \ref{s4}, which is based on \cite{Sa-St}, this is shown explicitly when $M$ is compact. In fact, some general causal and geometrical properties (Theorem \ref{t4a}, Proposition \ref{p4a}) are shown to imply the connectivity by  timelike geodesics (Theorem \ref{t4c}) and the existence of a closed timelike geodesic (Theorem \ref{t4d}); this  improves widely the analogous variational results in \cite{CMP}. In Section \ref{s5}, general properties of geodesics are discussed, including related problems as periodic trajectories or boundaries. Finally, in Section \ref{s6}, we show that the variational techniques applied in \cite{BCFS} (firstly introduced for the Riemannian case in \cite{CFS}) yield accurate results on geodesic connectedness, improving previous results and techniques.

\section{Completeness and standard static spacetimes} \label{s2}

In what follows, $(M,g)$ denotes a  connected $n(\geq 2)$-Lorentzian manifold $(-,+,$ $\cdots , +)$, which is smooth enough ($C^\infty$ as a simplification, but we only need $C^1$ for most purposes).  Our notation and conventions will be standard, as in the books \cite{BEE}, \cite{O}, \cite{SW}. A Lorentzian manifold will be called a  stationary spacetime if it admits a timelike Killing vector field $K$, and {\em static} if, additionally, $K$ is  irrotational (the orthogonal distribution to $K$ is involutive). In this case, we will refer to $K$ as the  static vector field, which will be assumed to time-orient the spacetime. $K$ is not determined univocally; nevertheless, if another static $K'$ exists and it satisfies
$K'= f K$ for some $f\in C^\infty (M)$, then $f$ is a constant \cite[p. 224]{SW}. Some properties of Lorentzian manifolds admitting a Killing vector field are reviewed in \cite{Sa-grec}.
A {\em standard} static spacetime is a product 
$\R\times S$ endowed with the metric 
\be \label{ems}
 g[(t,x)] = -\beta(x) dt^2 + g_S[x] 
\ee
where $g_S$ is a Riemannian metric on $S$. Locally, any static spacetime looks like a standard one, with $K$ identifiable to $\partial_t$. Note that function $\beta = -g(K,K) \in C^\infty(M)$  makes sense also in the non-standard case (as well as the integral hypersurfaces of the orthogonal foliation $K^\perp$). 

An auxiliary Riemannian metric $g_R$ can be constructed by reversing
the sign of $g(K,K)$, that is, 
$ g_R(u,v) = g(u,v)+ 2 g(u,K)g(v,K)/\beta , $ for all $u,v \in TM. $
The qualitative behavior of $\beta$ at infinity for this Riemannian metric (or, eventually, $g_S$ in the standard case) will be especially relevant. A function $V$ on a complete Riemannian manifold behaves {\em at most quadratically} (at infinity) if, for some $A >0$:
\be \label{e2} V(x) \leq A\cdot d(x)^2 
\ee
for all $x$ outside a compact subset, where $d$ is the distance (canonically associated to the Riemannian metric) to some fixed point $x_0$. 

The (geodesic) completeness of a static spacetime is interesting not only in its own right, but also in relation to the structure of the spacetime. 

\begin{teor} \label{t2a} Let $(M,g)$ be a static spacetime:

(1) If $g$ is complete, then both the static vector field $K$ and the integral hypersurfaces of 
$K^\perp$ are complete with the restriction of the metric $g$.

(2) If the auxiliary Riemannian metric $g_R$ is complete and $\beta^{-1}$ behaves at most quadratically,  then $g$ is complete. In particular, the following are complete: (a) compact static spacetimes, and (b) standard static spacetimes with $g_S$  complete and $\beta^{-1}$ at most quadratic on $(S,g_S)$. 

\end{teor}
{\em Sketch of proof.} For {\it (1)}, see \cite[Theorem 2.1(2)]{Sa-St}. For {\it (2)}, following \cite[Proposition 2.1]{RS2}, let $\gamma: [0,1)\rightarrow M$ be a $g$-geodesic. It is enough to prove that the $g_R$-length of $\gamma$ is bounded, (see \cite[Lemma 5.8]{O}). As $K$ is Killing, $g(\dot \gamma,K)= \lambda$ (constant) and $|\dot \gamma|^2_R (:=g_R(\dot \gamma,\dot \gamma))= C+2\lambda^2/\beta\circ \gamma $, with $C=g(\dot \gamma,\dot \gamma)$. Thus, choosing $x_0=\gamma(0)$ in the definition of $d$ below (\ref{e2}), as well as appropiate constants $B,D>0$:
$$ |\dot \gamma|_R(s) \leq B\cdot d(\gamma(s)) + D \leq B\cdot  \int_0^s 
|\dot \gamma|_R (\bar s)d\bar s + D , $$ 
 and, thus,
$$ \log \left(B\cdot  \int_0^s 
|\dot \gamma|_R (\bar s)d\bar s + D \right) - \log D \leq  B s \leq B .$$
Therefore, the length of $\gamma$ is bounded, as required. 
Finally, {\it (a)} is obvious and, for {\it (b)}, recall \cite[Proposition 7.40]{O}. $\quad \quad $
\cvd

\noindent {\bf Note.} (1) The proof of the first assertion in {\it (2)} holds in the stationary case (and {\it (a)} holds even for conformally stationary spacetimes). (2) The optimality of the at most quadratic behaviour for $\beta^{-1}$ can be explicitly checked in the standard case with $S=\R$.

\smallskip

\noindent Completeness of $K$ is an obvious necessary condition to ensure that a static spacetime is standard with $K=\partial_t$. Notice that this condition is very mild and automatically satisfied when $M$ is compact.
In the simply connected case, it is also sufficient \cite[Theorem 2.1(1)]{Sa-St}:

\begin{teor} \label{t2b}
A simply connected static spacetime is standard if and only if at least one of its static vector fields $K$ is complete.
\end{teor}

\noindent {\bf Note.} No analogous result holds in the stationary case. In fact, there exist simply connected compact stationary manifolds of any dimension $n\geq 5$ and $n=3$ (the exceptions $n=2,4$ appear only by topological reasons: no Lorentzian manifolds of such dimensions can be compact and simply connected).

\smallskip

\noindent Then, as a consequence of Theorems \ref{t2a}, \ref{t2b} we have (for further study in the non-simply connected case, see \cite{GO}):

\begin{coro} \label{c2a}
The universal covering of any geodesically complete static spacetime is a standard static spacetime.
\end{coro}

\section{Causality} \label{s3}

Causality of static spacetimes is simplified by considering the conformal metric 
$g^* = g/\beta$. In the standard case, $t$ is a time function (which implies causal stability) and $g^*$ is globally a product metric. Moreover:

\begin{teor} \label{t3a}
Any standard static spacetime $(M,g)$ as in {\rm (\ref{ems}) }
is causally continuous, and they are equivalent:

(A) $(M,g)$ is globally hyperbolic.

(B) The conformal metric $g_S^* = g_S/\beta$ is complete.

(C) Each slice $t=$constant is a Cauchy hypersurface.

In particular, this holds if $g_S$ is complete and $\beta$ at most quadratic.
\end{teor}
{\em Sketch of proof.} For the first assertion,  by \cite[Th. 3.25, Prop. 3.21]{BEE}, it is enough to prove past and future reflectivity. Let us show past reflectivity  $I^+(p) \supseteq I^+(q) \Rightarrow I^-(p) \subseteq I^-(q)$ (future is analogous). Put $p=(t_p,x_p), q=(t_q,x_q)$. Assuming the first inclusion, it is enough to prove $p_{-\epsilon} := (t_p-\epsilon,x_p) \in I^-(q)$, for all $\epsilon >0$. As $q_\epsilon := (t_q + \epsilon , x_q) \in I^+(p)$, there exists a future-directed timelike curve $\gamma(s)=(s,x(s)), s\in [t_p,t_q + \epsilon]$ joining $p$ and $q_\epsilon$. Then, the future-directed timelike curve $\gamma_{-\epsilon}(s)=(s-\epsilon  ,x(s))$ connects $p_{-\epsilon}$ and $q$, as required.

For the remainder, use   \cite[Theorems 3.67, 3.69]{BEE}. $\triangle$ 

\begin{remas}\label{r3a} {\rm
(1) When $g_S^*$ is not complete, the spacetime may be non-causally simple. In fact,  Take  $(S,g^*_S) = \R^2\backslash \{(1,y) : y \leq 1\}$, and choose $p=(0,0,0), q=(\sqrt 8, 2, 2) \in \R \times S$. Clearly 
$q\in  \overline{J^+(p)} \backslash J^+(p)$, that is, $J^+(p)$ is not closed  as required.  

(2) The required behaviour at infinity for $\beta$ is different in Theorems \ref{t3a} and \ref{t2a}, which shows the independence between completeness and global hyperbolicity. More striking, even if $g_S$ is incomplete (and, thus, so is $g$ for any $\beta$) a conformal factor $\Omega $ exists such that  $g_S^* = \Omega g_S$ is complete \cite{NO}. Thus, the corresponding $g=-\beta dt^2+g_S$  will be globally hyperbolic choosing $\beta = 1/\Omega$.

}\end{remas}
\noindent For further results  on causality in the standard case, see \cite{Sa-ba}.

In the non-standard one, recall first that a stationary spacetime is called {\em standard} when its metric can be written as in (\ref{ems}) plus, eventually, cross terms between $\R$ and $S$ independent of $t$. 

\begin{propo}\label{t3c}
A globally hyperbolic spacetime $(M,g)$ which admits a complete stationary vector field $K$ is standard stationary.
\end{propo}
{\em Proof.} By using \cite{BS}, a smooth spacelike Cauchy hypersurface $S$ can be found. Then, recall that $\R \times S \rightarrow M, (t,x) \rightarrow \Phi_t(x)$ is an isometry, where $\Phi$ is the flow of $K$. $\quad $
\cvd

\begin{remas} \label{r3b} {\em (1) In the proof each hypersurface $t=$ constant is  Cauchy and, again, the completeness of $K$ is essential (consider an open diamond $I^+(p) \cap I^-(q)$ in Lorentz-Minkowski space).

(2)  Proposition  \ref{t3c} holds in particular in the static case, but:
{\em in general,  
a complete static globally hyperbolic spacetime $(M,g)$ is not a standard static one.}  
To construct an explicit counterexample, take the cylinder $\R\times S^1$,  
$g= -E(\theta )dt^2 + 2 F(\theta ) dt d\theta + G(\theta) d\theta^2$, $E,  G >0$, which is complete and admits the Cauchy hypersurface $t=0$. Clearly,  $K=\partial_t$ is a static vector field (the distribution $K^\perp$ is one-dimensional), and it is unique (up to a constant) 
for choices of  $E, F, G$ with non-constant curvature (no two independent Killing vector fields can exist on a neighborhood with non-constant curvature, 
see the proof of \cite[Theorem 4.2(1)]{Sa-trans}). But, for most choices of $E,F,G$, the integral curves of $K^\perp$ are not circles and, thus, neither Cauchy hypersurfaces (see \cite[Example 2.4]{Sa-St} for related computations). 
}\end{remas}

\section{Compact static spacetimes} \label{s4}

In the bidimensional case, compact static spacetimes are essentially globally conformally flat Lorentzian tori, exhaustively studied in \cite{Sa-trans}. For the higher dimensional case, causality plays an important role.
 
\begin{teor}\label{t4a} Let $(M,g)$ be a compact static spacetime. Then:

(1) Its universal Lorentzian covering $(\bar M,\bar g)$ is standard static and admits a (possibly non-compact) Cauchy hypersurface.

(2) $(M,g)$ is totally vicious, i.e., $I^+(p)=I^-(p)=M, \forall p \in M$.

\end{teor}
{\em Sketch of proof.} {\it (1)} Straightforward from Theorem \ref{t2a}(2), Corollary \ref{c2a}, Theorem \ref{t3a} and the example in Remark \ref{r3b}(2). {\it (2)} The following general result (which includes the stationary case) holds \cite[Theorem 1.1]{Sa-St}:
{\em any compact Lorentzian manifold $(M,g)$ which admits a timelike conformal vector field  is totally vicious}. $\quad \quad $
\cvd

\noindent Even more,  deck transformations of $(\bar M,\bar g)$ satisfy:
\begin{propo} \label{p4a}
Let $(M,g)$ be a compact static manifold, and $(\bar M,\bar g)$ its universal Lorentzian covering. Then any deck transformation $\phi: \bar M \rightarrow \bar M$ can be written as
$$
\phi(t,x)= (t+T_\phi , \phi^S(x)),
$$
for some diffeomorphism  $\phi^S$ of $S$ and $T_\phi \in \R$.

\end{propo}

The striking difference between the causality of $M$ and its universal covering $\bar M$, yields interesting consequences for geodesic connectedness and existence of closed geodesics. Concretely, for the first question, one can ensure not only geodesic connectedness but also the existence of  timelike connecting geodesics:
\begin{teor} \label{t4c}
Any pair of points $p,q $ in a compact static spacetime $(M,g)$
 can be joined by means of a timelike geodesic.
\end{teor}
{\em Sketch of proof.} Join $p, q$ by means of a timelike curve $\gamma$ (Theorem \ref{t4a}(2)), and lift this curve to $\bar \gamma$ in the universal covering. Because of global hyperbolicity (Theorem \ref{t4a}(1)), its extremes $\bar p, \bar q$ can be joined by a timelike geodesic (Avez-Seifert theorem), which projects onto the required one. $\quad \quad$ 
\cvd

\noindent For the second question, it is known (Tipler \cite{Ti}, see also \cite[Theorem 4.15]{BEE}): {\em any compact spacetime regularly covered by a spacetime which admits a {\em compact} Cauchy hypersurface, contains a closed timelike geodesic.} An extension of Tipler's technique yields:

\begin{propo} \label{p4b}
Let $(M,g)$ be a compact Lorentzian manifold which admits a regular Lorentzian covering $\Pi: \bar M\rightarrow M$ such that $\bar M$ is globally hyperbolic. Assume that a 
conjugacy class ${\cal C}$ of the group of deck transformations of $\bar M$ is finite 
and contains a closed timelike curve $\alpha$.
Then there exists at least one closed timelike geodesic in ${\cal C}$. 
\end{propo}
A particular case of this result is \cite[Corollary 4.7]{CMP}.
Even more, Proposition \ref{p4b} admits further extensions, which are applicable whenever the deck transformations satisfy the conclusion of Proposition \ref{p4a}. Then \cite{Sa-St}:

\begin{teor} \label{t4d}
For a compact static spacetime $(M,g)$,  
any conjugacy class ${\cal C}$ of the group of deck transformations, containing a closed timelike curve,  admits a closed timelike geodesic. In particular, there exists at least one closed timelike geodesic in $M$.
\end{teor}

\section{Geodesic equations} \label{s5}

For simplicity we will consider in this section the standard case $\R \times S$, but extensions to the non-standard one are obvious. From a straightforward computation, $\gamma(s) =(t(s),x(s))$ is a geodesic if and only if:
\be \label{e5a}
D_s \dot{ x} = -\frac{1}{2}\dot{t}^2 \nabla \beta (x) \, , \quad \quad \quad  \beta (x) \dot t = \lambda \, , 
\ee
for some $\lambda \in \R$.
When $\lambda =0$ then $\gamma(s) =(t_0, x(s))$, for some $t_0\in \R$, is essentially a geodesic on $S$. Otherwise, we can normalize the value of $\lambda$, and the following result holds (see \cite{Sa-nonlin} for details):

\begin{propo} \label{p5a}
Let $\gamma(s) =(t(s),x(s))$ be a  geodesic of $S$ satisfying {\rm (\ref{e5a})} with $\lambda = \sqrt{2}$ and $E=g(\dot \gamma,\dot \gamma)/2$.  Let $V=-1/\beta$, then:

(1) $x(s)$ is a solution of classical equation
$ D_s \dot{ x} = -\nabla V(x)$, with ``total energy'' $E=g_S(\dot x, \dot x)/2+ V(x)$.

(2) $x(s)$ is a pregeodesic (i.e., geodesic  up to a reparameterization) for the Jacobi metric 
$g_E = (E - V) g_S$, whenever $E>V$.
\end{propo}
Conversely, from any curve satisfying {\it (1)} or {\it (2)}, a solution to (\ref{e5a}) with $\lambda= \sqrt{2}$ (i.e., a normalized geodesic) can be constructed.
This yields an equivalence between the set of geodesics (non-orthogonal to $K$) 
in the static spacetime and the trajectories for the potential $V$, which can be applied to some problems \cite{Sa-nonlin}. In particular, to look for geodesics $\gamma(s)=(t(s),x(s))$ with ``energy'' (rest mass, if non-spacelike) $-2E=-g(\dot \gamma,\dot \gamma)$, which are ``spacelike closed'' in the following sense: for some $b>0$ (proper period) and $ T>0$ (universal period), $ x(b)=x(0), \dot x(b)=\dot x(0), t(b)=t(0)+T$. Proposition \ref{p5a} applies directly to find such geodesics with fixed value of $E$ or $b$; for fixed $T$ (``T-periodic trajectories''), see \cite{Sa-pams} and references therein.

Finally, we point out that, starting with Proposition \ref{p5a}, the equivalence between two  notions of {\em convexity} for incomplete $g_S$ can be proven \cite{BGS}. When $(S,g_S)$ can be regarded as an open subset with smooth boundary $\partial S$ of a bigger complete Riemannian manifold, there are some natural Riemannian notions of convexity, essentially equivalent. Nevertheless, their equivalence for a Lorentzian manifold, as a static spacetime $(\R \times S, g)$, has to be proven (see \cite{Sa-cata} for a review). Even more, convexity for each causal character (timelike, lightlike, spacelike) must be taken into account. In the static case, the problem is solved in \cite{BGS}, proving in particular the equivalence between the geometrical and  variational notions of convexity, which are useful for the problem of geodesic connectedness (see also \cite{BCFS} for non-smooth boundaries).   

\section{Variational approach. Geodesic connectedness} \label{s6}

A different approach is required for geodesic connectedness  (see \cite{Sa-cata} for a survey on this problem). This  will be discussed only in the standard case $M=\R\times S$ because, essentially, the natural hypotheses for connectedness are applicable to the universal covering, and Theorem \ref{t2b} would be also applicable. 
A fruitful idea is the variational characterization of connecting geodesics by Benci, Fortunato and Giannoni \cite{BFG} (see also \cite{Ma}):  given two points $p_0=(t_0,x_0), p_1=(t_1,x_1)\in \R\times S$, a (absolutely continuous) curve $\gamma: [0,1] \rightarrow \R\times S, $ $\gamma(s)=(t(s),x(s))$ with $p_0=\gamma(0), p_1=\gamma(1)$ is a geodesic if and only if the $x(s)$ component is a critical point of the functional:
\begin{equation}
\label{e6a}
J(x) =\ {1\over 2} \int_0^1 g_S(\dot x(s),\dot x(s)) ds\
 -\ \frac{\Delta_t^2}{2}\ \left(\int_0^1 {1\over \beta(x(s))}\ ds\right)^{-1},
\end{equation}
where $\Delta_t^2 = (t_1-t_0)^2$. Thus, the problem is reduced to consider the functional $J$ on curves $x:[0,1]\rightarrow S$ connecting $x_0, x_1$. Even more, these authors studied directly the functional $J$, ensuring the existence of critical points if: (i) $g_S$ is complete, and (ii)  $\beta$ is (upper) bounded. Recently, Caponio, Masiello and Piccione \cite{CMP} showed that  assumption (ii) can be weakened in: (ii') $\beta$ is subquadratic (i.e., inequality (\ref{e2}) holds replacing $d^2$ by $d^{2-\epsilon}$, for some $\epsilon >0$). On the other hand, Allison and Unal \cite{AU} showed that, when (i) and (ii) holds, if $S$ is pseudoconvex and non-returning, then so is $M$; in this particular case, an alternative non-variational proof of geodesic connectedness is obtained.

Recently, we \cite{BCFS} have obtained a fully satisfactory approach to this problem, which can be summarized in the following points:
\begin{enumerate}
\item By using Cauchy-Schwarz inequality, $J$ is lower bounded by
\begin{equation}\label{e6b}
J(x)  \geq  \frac{1}{2}\int_0^1 g_S(\dot x(s),\dot x(s))\ ds\
 -\ \frac{\Delta_t^2}{2}\ \int_0^1  \beta(x) \ ds ,
\end{equation}
which can be seen as the action $A= \int L$ associated to a typical Lagrangian type $L=T-V$, with $T=g_S(\dot x,\dot x)/2$ (kinetic energy) and $V= \Delta_t^2 \beta /2$ (potential energy). General properties which imply the existence of critical points for a functional (coercitivity, boundedness from below), if are fulfilled by $A$, then are also fulfilled by $J$.

\item These general properties for action $A$ were known to hold when $g_S$ is complete and $V$ subquadratic (this re-proves Caponio et al. result). Even more, an optimal result for $V$ at most quadratic was obtained in \cite{CFS}. Nevertheless, this result does not guarantees the existence of critical points for $V$ quadratic (in fact, there are simple counterexamples) and, thus, neither the geodesic connectedness of $M$.

\item Studying directly the functional $J$ for $\beta$ quadratic (without the loss of information in Cauchy-Schwarz inequality), by using the accurate technique in \cite{CFS},  both coercitivity and lower boundedness for $J$, are ensured, i.e.:

\end{enumerate}

\begin{teor} \label{t6a}
A standard static spacetime with $g_S$ complete and $\beta$ at most quadratic
 is geodesically connected.
\end{teor}

\begin{remas} {\rm (1)
The assumption on completeness for $g_S$ is unavoidable (a strip of the
universal covering of two--dimensional anti--de Sitter spacetime
such as $M= \R \times (-\pi/4, \pi/4) $, $g = (dx^2-dt^2)/\cos^2x$ is not geodesically connected). Nevertheless, completeness can be weakened by the existence of a con\-vex boundary. 
The ``at most quadratic'' behaviour for $\beta$ is optimal, as shown by a family of explicit counterexamples \cite[Section 7]{BCFS} (inspired by \cite[Section 6.1]{FS}).

(2) Causally  related points will be connectable by (maximizing) causal geodesics because the assumptions imply global hyperbolicity (Theorem \ref{t3a}). But causal connectability will also hold if the completeness of $g_S$ is weakened by the time and lightlike convexities of the boundary. 

(3) A related variational problem is the connectability by means of trajectories under a potential vector (as charged particles), see for example, \cite{Min}. In particular, when applied to static spacetimes, this also yields implications on connectability for stationary ones.
}
\end{remas}



\begin{thebibliography}{00}



\bibitem{AU} \sc D.E. Allison, B. Unal,  {\em Geodesic structure of standard static spacetimes}, {\bf 863}, \rm J. Geom Phys., {\bf 46} 193--200 (2003).



\bibitem{BCFS} \sc R. Bartolo, A. M. Candela, J.L. Flores, M. S\'anchez, 
{\em Geodesics in static Lorentzian manifolds
with critical quadratic behavior}, \rm  Adv. Nonlinear Stud. {\bf 3}  471-494 (2003).


\bibitem{BGS}
\sc R. Bartolo, A. Germinario, M. S\'anchez, {\it A note  on the boundary of  a static Lorentzian
manifold}, \rm Differential Geom. Appl. {\bf 16}  121-131 (2002).



\bibitem{BEE} \sc J.K. Beem, P.E. Ehrlich and  K. Easley, \rm Global  Lorentzian  Geometry, Marcel Dekker Inc., N.Y. (1996).



\bibitem{BFG}
\sc V. Benci, D. Fortunato and F. Giannoni,
\em On the existence of multiple geodesics in static space--times,
\rm Ann. Inst. H. Poincar\'e Anal. Non Lin\'eaire 
{\bf 8} 79-102 (1991).

\bibitem{BS} \sc A. N. Bernal, M. S\'anchez, \rm {\em On Smooth Cauchy Hypersurfaces and Geroch's Splitting Theorem},    Commun. Math. Phys., {\bf 243}  461-470 (2003). 


\bibitem{CFS} \sc A. M. Candela, J.L. Flores, M. S\'anchez, {\em A Bolza-type problem in a Riemannian manifold}, \rm J. Differential Equations, {\bf 193} 196-211 (2003).




\bibitem{CMP} \sc E. Caponio, A. Masiello, P. Piccione,  {\em Some global properties of static spacetimes}, 
\rm Math. Z., {\bf 244}, 457--468  (2003).



\bibitem{FS} \sc J.L. Flores, M.  S\'anchez, {\em Geodesic connectedness and conjugate points in GRW space-times}, \rm J. Geom. Phys. {\bf 36} 285--314  (2000). 


\bibitem{GO} \sc M. Guti\'errez, B. Olea, \em Splitting theorems in the presence of an irrotational vector field, \rm Univ. M\'alaga, preprint.

\bibitem{Ma} 
\sc A. Masiello,
\rm Variational methods in Lorentzian geometry,
\rm Pitman Res. Notes Math. Ser. {\bf 309},
Longman Sci. Tech., Harlow, 1994.

\bibitem{Min} \sc E. Minguzzi, \em On the existence of maximizing curves for the charged-particle action, \rm Classical Quantum Gravity {\bf 20} 4169--4175, (2003).

\bibitem{NO} \sc K. Nomizu, H. Ozeki, \em The existence of complete Riemannian metrics, \rm Proc. Amer. Math. Soc. {\bf 12} 889-891 (1961).  


\bibitem{O} \sc B. O'Neill, \rm Semi-Riemannian Geometry, Ser.  Pur.  Appl. Math. {\bf 103}, 
Academic Press, New York, (1983).

\bibitem{RS2}  \sc A. Romero, M. S\'anchez, {\em On the completeness of certain families  of semi-Riemannian manifolds}, \rm Geometriae Dedicata {\bf 53}, 103-117 (1994).

\bibitem{SW} \sc R. Sachs, H. Wu,  \rm
General Relativity for Mathematicians, Graduate Texts in Math. {\bf 48}, Springer, New York (1977).

\bibitem{Sa-St} \sc M. S\'anchez, 
{\em On causality and closed geodesics of compact Lorentzian manifolds and static spacetimes}, \rm Univ. Granada, preprint.

\bibitem{Sa-cata} \sc M. S\'anchez, 
{\em Geodesic connectedness of semi-Riemannian manifolds},
\rm Nonlinear Anal., {\bf 47},  3085-3102, (2001).

\bibitem  {Sa-nonlin} \sc M.  S\'anchez, \rm {\em Geodesics  in  static  spacetimes   and   t-periodic trajectories}
\rm Nonlinear Anal.,  {\bf 35}, 677-686 (1999).

\bibitem  {Sa-pams}  \sc M. S\'anchez, {\em Timelike periodic trajectories in spatially compact Lorentzian manifolds}, \rm Proc.  Amer.  Math.  Soc. {\bf 127},  3057-3066 (1999).

\bibitem  {Sa-trans} \sc M.  S\'anchez, 
{\em Structure of Lorentzian tori with a Killing vector field}, 
\rm Trans. Amer. Math. Soc., {\bf 349}, 1063-1080 (1997). 

\bibitem{Sa-grec} \sc M. S\'{a}nchez, {\em  Lorentzian Manifolds Admitting a Killing Vector 
Field},  \rm Nonlinear Anal., {\bf 30}  643-654 (1997). 

\bibitem{Sa-ba} \sc M. S\'anchez, {\em
	Some remarks on Causality and Variational Methods in Lorentzian 	manifolds, }	
    {\rm  Conf. Sem. Mat. Univ. Bari} {\bf 265} 1-12 (1997).

\bibitem{Ti} \sc F.J. Tipler, {\em Existence of a closed timelike geodesic in Lorentz spaces}, {\rm Proc. Amer. Mah. Soc.} {\bf 76,} 145-147 (1979). 







\end{thebibliography}
\end{document}